\documentclass{article} 

 \usepackage[usenames,dvipsnames]{pstricks}
 \usepackage{epsfig}
 \usepackage{pst-grad} 
 \usepackage{pst-plot} 
 \usepackage[space]{grffile} 
 \usepackage{etoolbox} 
 \makeatletter 
 \patchcmd\Gread@eps{\@inputcheck#1 }{\@inputcheck"#1"\relax}{}{}
 \makeatother

\usepackage[usenames,dvipsnames]{pstricks}
\usepackage{epsfig}
\usepackage{pst-grad} 
\usepackage{pst-plot} 
\usepackage[space]{grffile} 
\usepackage{etoolbox} 
\makeatletter 
\patchcmd\Gread@eps{\@inputcheck#1 }{\@inputcheck"#1"\relax}{}{}
\makeatother

 \usepackage[usenames,dvipsnames]{pstricks}
 \usepackage{epsfig}
 \usepackage{pst-grad} 
 \usepackage{pst-plot} 
 \usepackage[space]{grffile} 
 \usepackage{etoolbox} 
 \makeatletter 
 \patchcmd\Gread@eps{\@inputcheck#1 }{\@inputcheck"#1"\relax}{}{}
 \makeatother

\usepackage{amsmath}
\usepackage{amssymb}
\usepackage{amsthm}
\usepackage[latin1]{inputenc}
     
\usepackage{graphicx}

\usepackage[usenames,dvipsnames]{pstricks}
\usepackage{epsfig}
\usepackage{pst-grad} 
\usepackage{pst-plot} 

\usepackage{ifthen}
\usepackage{color}

\newcommand{\intav}[1]{\mathchoice {\mathop{\vrule width 6pt height 3 pt depth  -2.5pt
\kern -8pt \intop}\nolimits_{\kern -6pt#1}} {\mathop{\vrule width
5pt height 3  pt depth -2.6pt \kern -6pt \intop}\nolimits_{#1}}
{\mathop{\vrule width 5pt height 3 pt depth -2.6pt \kern -6pt
\intop}\nolimits_{#1}} {\mathop{\vrule width 5pt height 3 pt depth
-2.6pt \kern -6pt \intop}\nolimits_{#1}}}

\def\polhk#1{\setbox0=\hbox{#1}{\ooalign{\hidewidth\lower1.5ex\hbox{`}\hidewidth\crcr\unhbox0}}}

\newcommand{\bmo}{\operatorname{BMO}}
\newcommand{\dist}{\operatorname{dist}}

\newcommand{\Id}{\operatorname{Id}}

\newcommand{\tr}{\operatorname{Tr}}

\newcommand{\llip}{\operatorname{Log-Lip}}

\newcommand{\bb}{{\bf b}}

\newtheorem{teo}{Theorem}[section]

\newtheorem{Definition}{Definition}[section]
\newtheorem{Lemma}{Lemma}[section]

\newtheorem{Proposition}{Proposition}[section]
\newtheorem{Remark}{Remark}[section]
\newtheorem{Assumption}{A}

\begin{document}

\title{Regularity theory for the Isaacs equation through approximation methods}
\author{Edgard A. Pimentel}

\date{\today} 

\maketitle

\begin{abstract}

\noindent In this paper, we propose an approximation method to study the regularity of solutions to the Isaacs equation. This class of problems plays a paramount role in the regularity theory for fully nonlinear elliptic equations. First, it is a model-problem of a non-convex operator. In addition, the usual mechanisms to access regularity of solutions fall short in addressing these equations. We approximate an Isaacs equation by a Bellman one, and make assumptions on the latter to recover information for the former. Our techniques produce results in Sobolev and H\"older spaces; we also examine a few consequences of our main findings.
\medskip


\medskip 

\end{abstract}

\vspace{.1in}

\section{Introduction}\label{introduction}

In the present paper, we examine the regularity of the solutions to an Isaacs equation of the form

\begin{equation}\label{eq_isaacs}
	\sup_{\alpha\in\mathcal{A}}\inf_{\beta\in\mathcal{B}}\,\left[-\tr\left(A_{\alpha,\beta}(x)D^2u\right)\right]\,=\,f\;\;\;\;\;\mbox{in}\;\;\;\;\;B_1,
\end{equation}
where $A_{\alpha,\beta}:B_1\times\mathcal{A}\times\mathcal{B}\to\mathbb{R}^{d^2}$ is a $(\lambda,\Lambda)$-elliptic matrix, $f\in \bmo(B_1)$ and $\mathcal{A}$ and $\mathcal{B}$ are compact, separable and countable metric spaces. 

We argue by an approximation method, relating \eqref{eq_isaacs} to a Bellman equation. As a result, we study the regularity of solutions to \eqref{eq_isaacs} in Sobolev and H\"older spaces. Further, we investigate a few consequences of our findings. 

The primary motivation for the study of the Isaacs equations comes from the two-players, zero-sum, (stochastic) differential games. See \cite{isaacsbook} and \cite{friedmanbook}. In this context, solutions to \eqref{eq_isaacs} are, at least formally, value functions for the associated differential game. This fact is consequential on the Dynamic Programming Principle, together with further considerations from game theory. Applications of this theory are in problems of competitive advertising, duopolistic competition and models of resource extraction; see, for example, \cite{yeung}.

The theory of viscosity solutions led to important advances regarding \eqref{eq_isaacs}. We refer the reader to \cite{crandalllions}, \cite{crandallevanslions} and \cite{lionsbook82}.  We emphasize the well-posedness for \eqref{eq_isaacs}, representation formulas, optimality conditions for the solutions and the existence of value functions for the associated game. See \cite{evanssouga1}, \cite{flemingsouganidis1}, \cite{swiech96}, \cite{buckdahn2} and \cite{buckdahn3} and the references therein. Representation formulas for viscosity solutions to (degenerate) parabolic equations is the subject of \cite{katsoulakis}. For a survey on the topic, we refer the reader to \cite{buckdahn1}. 

The interest in the regularity of the solutions to \eqref{eq_isaacs} is in two aspects. First, it is an example of a non-convex/non-concave equation; hence, the Evans-Krylov theory does not apply to \eqref{eq_isaacs}. In fact, the best results holding in general are estimates in $\mathcal{C}^{1,\gamma}$ and $W^{2,\mu}$, for some $\gamma$ and $\mu$, universal. We refer the reader to \cite{krysaf1}, \cite{krysaf2} and \cite{fanghuaw2delta}. See further \cite{ccbook}. In addition, any fully nonlinear elliptic equation can be written as an Isaacs operator; see \cite[Remark 1.5]{cc2003}. This fact stems from the ellipticity and the structure of the Pucci's extremal operators. 

The regularity theory for the equations of the form \eqref{eq_isaacs} is launched in \cite{cc2003}. In that paper, the authors study
\begin{equation}\label{eq_cabcaf1}
	F(D^2u)\,:=\,\min\left\lbrace F_\wedge(D^2u),F_\vee(D^2u) \right\rbrace\,=\,0\;\;\;\;\;\mbox{in}\;\;\;\;\;B_1,
\end{equation}
where $F_\wedge$ and $F_\vee$ are, respectively, concave and convex elliptic operators. In this case, they prove that solutions are in $\mathcal{C}^{2,\gamma}_{loc}(B_1)$. In addition, the authors establish the existence of solutions to the Dirichlet problem for \eqref{eq_cabcaf1} with merely continuous boundary data. 

A similar problem is examined in \cite{collins}. In that paper, the author produces a priori estimates for the solutions to
\begin{equation}\label{eq_twisted}
	F(D^2u)\,:=\, F_\wedge(D^2u)\,+\,F_\vee(D^2u)\,=\,0\;\;\;\;\;\mbox{in}\;\;\;\;\;B_1.
\end{equation}
In particular, the author shows that classical solutions have estimates in $\mathcal{C}^{2,\gamma}$.

As regards nonconvex equations, we mention \cite{caffyu}. In that paper, the authors consider an homogeneous fully nonlinear equation governed by a nonconvex operator. Under assumptions on the level sets of the operator the authors prove that smooth solutions have estimates in $\mathcal{C}^{2,\gamma}$.

The regularity of the solutions to \eqref{eq_isaacs} in Sobolev spaces is the subject of \cite{kovats12}. The author supposes that $A_{\alpha,\beta}$ is a diagonal, separable matrix, Lipschitz-continuous with respect to $\alpha$ and $\beta$ and independent of $x$. In this setting, he proves that classical solutions are in $W^{2,p}(B_1)$ for every $p>1$. Since solutions are required to be of class $\mathcal{C}^2$, this result is closer to the realm of \textit{a priori estimates} than to the regularity theory. We refer the reader to \cite{silvsirak} for a neat remark on the distinction between \textit{a priori estimates} and \textit{regularity theory}. 

H\"older continuity for the gradient of the solutions is the subject of \cite{krylov14}. In that paper, the author proves that solutions are in $\mathcal{C}^{1,\gamma}_{loc}(B_1)$, for some universal $\gamma\in(0,1)$. In addition, the author proves a convergence rate for a finite-differences scheme that approximates solutions. 

As regards approximating solutions to \eqref{eq_isaacs}, we mention \cite{cafsilsmo}. In that paper, the authors define a family $(E_\varepsilon)_{\varepsilon\in\mathcal{E}}$ of approximate operators; these are inspired by the theory of nonlocal equations. The authors prove that solutions $u_\varepsilon$ to the PDE governed by $E_\varepsilon$ are in $\mathcal{C}^{1,\gamma}_{loc}(B_1)$, uniformly in $\varepsilon>0$. Moreover, $u_\varepsilon$ is a classical solution for every $\varepsilon>0$. Finally, they obtain a convergence rate of the form
\[
	\left\|u_\varepsilon\,-\,u\right\|_{L^\infty(B_1)}\,\leq\,C\varepsilon^\nu,
\]
for some small $\nu>0$, where $u$ solves \eqref{eq_isaacs}.

In \cite{kovats16}, the author examines a very particicular example of \eqref{eq_isaacs}. Under conditions on the geometry of the domain, he proves estimates in H\"older and Sobolev spaces.

In spite of significant developments, the regularity theory for the Isaacs equation is far from being complete. In \cite{nadvla11}, the authors produce an example of a singular solution to the Isaacs equation. Then, estimates in $\mathcal{C}^{1,1}$ are not available in general.

In \cite{armstrong12}, the authors consider fully nonlinear $(\lambda,\Lambda)$-elliptic operators $F$ which are differentiable at the origin. Under this assumption, they prove the existence of a number $\varepsilon=\varepsilon(\lambda,\Lambda,d)$ and a set $\Omega\subset \overline{B}_1$ such that: \textit{i.} solutions to $F(D^2u)=0$ satisfy $u\in\mathcal{C}^{2,\gamma}(B_1\setminus\Omega)$ and \textit{ii.} the Hausdorff dimension of $\Omega$ is at most $d-\varepsilon$. Unfortunately, this (partial) result does not apply to the Isaacs equation, as mentioned by the authors. Being homogeneous of degree one, were the Isaacs operator differentiable at the origin, it would be linear.

An important aspect of \eqref{eq_isaacs} is that most mechanisms to access regularity fail in the context of this equation. For example, consider the geometric method introduced in \cite{caffarelli89}. In that paper, the author approximates a fully nonlinear elliptic operator $F(M,x)$ by its counterpart with fixed coefficients $F(M,x_0)$. Then, assumptions are imposed on the latter. The goal is to import regularity from $F(M,x_0)$ to the original problem. Typically, the assumptions placed on $F(M,x_0)$ are satisfied when this operator is convex. In the case of the Isaacs equation, this condition fails to hold for the fixed coefficients operator is of the form
\[
	\sup_{\alpha\in\mathcal{A}}\inf_{\beta\in\mathcal{B}}\,\left[-\tr\left(A_{\alpha,\beta}(x_0)M\right)\right].
\]

More recently, and in the spirit of \cite{caffarelli89}, an asymptotic technique was introduced. In \cite{silvtei}, the authors consider a fully nonlinear elliptic operator $F:\mathcal{S}(d)\to\mathbb{R}$ and define the \textit{recession function} associated with $F$ as follows:
\[
	F^*(M)\,:=\,\lim_{\mu\to 0}\mu F(\mu^{-1}M).
\]
In this approach, assumptions are placed on the recession operator $F^*$. In the case of the Isaacs equation, it falls short. Because this equation is positive homogeneous of degree one, the recession function recovers the Isaacs operator, yielding no further information. 

In the present paper, we approximate \eqref{eq_isaacs} by a Bellman equation of the form
\begin{equation}\label{eq_belint}
	\inf_{\beta\in\mathcal{B}}\,\left[-\tr\left(\overline{A}_\beta(x)D^2v\right)\right]\,=\,0\;\;\;\;\;\mbox{in}\;\;\;\;\;B_1,
\end{equation}
where $\overline{A}_\beta:B_1\times\mathcal{B}\to\mathbb{R}^{d^2}$ is a $(\lambda,\Lambda)$-elliptic matrix. We work under a smallness regime for the quantity
\[
	\left|A_{\alpha,\beta}(x)\,-\,\overline{A}_{\beta}(x)\right|;
\]
distinct smallness regimes yield different regularity results; see Section \ref{sec_assumptions}.

We observe that, for every $x_0\in B_1$, the Bellman equation \eqref{eq_belint} is convex with respect to the Hessian. Therefore, the Evans-Krylov theory is available and solutions are locally of class $\mathcal{C}^{2,\gamma}$, for some universal $\gamma\in(0,1)$. At the core of our arguments is the idea of importing regularity from a Bellman equation to the Isaacs one. 

Our first result regards estimates in Sobolev spaces. It includes operators with explicit dependence on the gradient. That is, we study equations of the form
\begin{equation}\label{eq_isaacslot}
	\sup_{\alpha\in\mathcal{A}}\inf_{\beta\in\mathcal{B}}\left[ -\tr\left(A_{\alpha,\beta}(x)D^2u\right)\,-\, \bb_{\alpha,\beta}(x)\cdot Du\right]\,=\,f\;\;\;\;\;\mbox{in}\;\;\;\;\;B_1,
\end{equation}
where $\bb_{\alpha,\beta}:B_1\times\mathcal{A}\times\mathcal{B}\to\mathbb{R}^d$ is a given vector field. We prove the following theorem:

\begin{teo}[Estimates in Sobolev spaces]\label{teo_main1}
Let $u\in\mathcal{C}(B_1)$ be a viscosity solution to \eqref{eq_isaacslot}. Let $d<p<q$. Suppose that Assumptions A\ref{assump_ellipticity}-A\ref{assump_vf}, to be detailed in Section \ref{sec_assumptions}, are in force. Then, $u\in W^{2,p}_{loc}(B_1)$ and
\[
	\left\|u\right\|_{W^{2,p}(B_{1/2})}\,\leq\,C\left(\left\|u\right\|_{L^\infty(B_1)}\,+\,\sup_{\alpha\in\mathcal{A}}\sup_{\beta\in\mathcal{B}}\left\|\bb_{\alpha,\beta}\right\|_{L^\infty(B_1)}\,+\,\left\|f\right\|_{L^p(B_1)}\right),
\]
where $C>0$ is a universal constant.
\end{teo}

The proof of Theorem \ref{teo_main1} unfolds along three main steps. First, our approximation methods build upon $W^{2,\mu}$-estimates. It produces $W^{2,p}$-regularity for the solutions of \eqref{eq_isaacs}, i.e., the equation without dependence on the gradient, for $d<p<q$. Then, we follow \cite{swiech97} and prove $W^{1,\overline{p}}$-estimates for the solutions of \eqref{eq_isaacslot}. Finally, a reduction argument closes the proof. As a consequence of Theorem \ref{teo_main1}, we produce estimates in John-Nirenberg spaces. See Remark \ref{rem_bmo}.

Our second result regards the borderline case. That is, the regularity of solutions in $\mathcal{C}^{1,\llip}_{loc}(B_1)$. 

\begin{teo}[Estimates in $\mathcal{C}^{1,\llip}$]\label{teo_main2}
Let $u\in\mathcal{C}(B_1)$ be a viscosity solution to \eqref{eq_isaacs}. Suppose that Assumptions A\ref{assump_ellipticity} and A\ref{assump_osc}, to be detailed in Section \ref{sec_assumptions}, are in force.
Let $x_0\in B_{1/2}$. Then, $u\in\mathcal{C}^{1,\llip}_{loc}(B_1)$ and
\[
	\sup_{x\in B_r(x_0)}\,|u(x)-[u(x_0)+Du(x_0)\cdot x]|\,\leq\, C\left(\left\|u\right\|_{L^\infty(B_1)}\,+\,\left\|f\right\|_{L^p(B_1)}\right)\left(-r^2\,\ln r\right),
\]
where $C>0$ is a universal constant and $0<r\leq 1/2$.
\end{teo}
As regards the borderline case, we use Theorem \ref{teo_main2} to produce estimates for the gradient in John-Nirenberg spaces. This result relies upon the notion of recession function and relates to \cite{cafsilsmo}.

Finally, we refine the smallness regime and examine estimates in $\mathcal{C}^{2,\gamma}$. In this setting, we access distinct regularity profiles by imposing further conditions on the source term. Our last theorem reads as follows:

\begin{teo}[Estimates in $\mathcal{C}^{2,\gamma}$]\label{teo_main3}
Let $u\in\mathcal{C}(B_1)$ be a viscosity solution to \eqref{eq_isaacs}. Suppose that Assumptions A\ref{assump_ellipticity} and A\ref{assump_oscc2alpha}, to be detailed in Section \ref{sec_assumptions}, are in force. 
\begin{enumerate}
\item Then, there exists $\gamma\in(0,1)$ such that $u$ is of class $\mathcal{C}^{2,\gamma}$ at the origin. 
\item If $f\equiv 0$, $u\in\mathcal{C}^{2,\gamma}_{loc}(B_1)$ and there is a universal constant $C>0$ such that
\[
	\left\|u\right\|_{\mathcal{C}^{2,\gamma}(B_{1/2})}\,\leq\,C\left\|u\right\|_{L^\infty(B_1)}.
\]
\end{enumerate}
\end{teo}

The remainder of this paper is structured as follows: in Section \ref{sec_elemnot} we collect preliminary definitions and elementary facts. Section \ref{sec_assumptions} details the main assumptions under which we work. We present the proof of Theorem \ref{teo_main1} in Section \ref{isaacsw2p}. The proof of Theorem \ref{teo_main2} is the subject of Section \ref{isaacsc1ll}, whereas Section \ref{isaacsc2alpha} presents the proof of Theorem \ref{teo_main3}. 

\bigskip

\noindent {\bf Acknowledgements:} We are grateful to Prof. Andrzej \'Swi\polhk ech for suggesting the use of approximation methods in the study of the Isaacs equation. We also thank Prof. Boyan Sirakov, Prof. Andrzej \'Swi\polhk{e}ch and Prof. Eduardo Teixeira for their interest, suggestions and valuable comments on this material. We are also grateful to an anonymous referee for his/her comments and suggestions, which improved substantially the material in this paper. We are partially supported by FAPERJ and PUC-Rio start-up and baseline funds.

\section{Preliminaries and main assumptions}

In what follows, we put forward some preliminary notions and material and detail our main assumptions.

\subsection{Elementary notions}\label{sec_elemnot}

Next, we collect a number of definitions and elementary results. Throughout the paper, $B_r$ denotes the ball of radius $r>0$ centered at the origin; $Q$ stands for the unit cube in $\mathbb{R}^d$. 

Let $0<\lambda\leq\Lambda$. An operator $F:\mathcal{S}(d)\to\mathbb{R}$ is said to be $(\lambda,\Lambda)$-elliptic if 
\[
	\lambda\|N\|\,\leq\,F(M\,+\,N)\,-\,F(M)\,\leq\,\Lambda\|N\|,
\]
holds for every $M,\,N\in\mathcal{S}(d)$ with $N\geq 0$. Analogously, a matrix $A:B_1\to\mathbb{R}^{d^2}$ is $(\lambda,\Lambda)$-elliptic if
\[
	\lambda \Id\,\leq\,A(x)\,\leq\,\Lambda \Id,
\]
for every $x\in B_1$. Next we introduce two definitions.

\begin{Definition}[Pucci's extremal operators]\label{def_pucci}
The Pucci's extremal operators $\mathcal{M}^\pm_{\lambda,\Lambda}:\mathcal{S}(d)\to\mathbb{R}$ are defined as follows:
\[
	\mathcal{M}^+_{\lambda,\Lambda}(M)\,:=\,\Lambda\sum_{e_i>0}e_i\,+\,\lambda\sum_{e_i<0}e_i\;\;\;\;\mbox{and}\;\;\;\;\mathcal{M}^-_{\lambda,\Lambda}(M)\,:=\,\Lambda\sum_{e_i<0}e_i\,+\,\lambda\sum_{e_i>0}e_i,
\]
where $e_i$ are the eigenvalues of $M$.
\end{Definition}

The extremal operators are of paramount importance for the theory of fully nonlinear elliptic equations. In particular, they allow us to define the class of $(\lambda,\Lambda)$-viscosity solutions.

\begin{Definition}[Class of viscosity solutions]\label{def_viscclass}
Let $f\in\mathcal{C}(B_1)$. A function $u\in\mathcal{C}(B_1)$ is in the class of supersolutions $\overline{S}(\lambda,\Lambda, f)$ if 
\[
	\mathcal{M}^-_{\lambda,\Lambda}(D^2u)\,\leq\, f\;\;\;\;\;\mbox{in}\;\;\;\;\;B_1,
\]
in the viscosity sense. In addition, $u\in\mathcal{C}(B_1)$ is in the class of subsolutions $\underline{S}(\lambda,\Lambda, f)$ if 
\[
	\mathcal{M}^+_{\lambda,\Lambda}(D^2u)\,\geq\, f\;\;\;\;\;\mbox{in}\;\;\;\;\;B_1,
\]
in the viscosity sense. The class of $(\lambda,\Lambda)$-viscosity solutions is the set 
\[
	S(\lambda,\Lambda,f)\,:=\,\overline{S}(\lambda,\Lambda, f)\cap\underline{S}(\lambda,\Lambda, f).
\]
\end{Definition}

Let $\mathcal{O}\subset\mathbb{R}^d$ be a $\mathcal{C}^{1,1}$ bounded domain. A paraboloid of opening $M\in\mathbb{R}$ is a function $P_M:\mathcal{O}\to\mathbb{R}$ of the form
\[
	P_M(x)\,:=\,\ell(x)\,+\,\frac{M}{2}|x|^2,
\]
where $\ell:\mathcal{O}\to\mathbb{R}$ is an affine function.

\begin{Definition}
Let $u\in\mathcal{C}(\mathcal{O})$. We set
\[
	\overline{\Theta}(u,\mathcal{O})(x)\,:=\,\inf_{M\in\mathbb{R}}\left\lbrace\exists \,P_M\,|\,P_M(x)=u(x)\;\;\mbox{and}\;\;P_M(y)\geq u(y),\;\;\forall y\in \mathcal{O}\right\rbrace
\]
and
\[
	\underline{\Theta}(u,\mathcal{O})(x)\,:=\,\inf_{M\in\mathbb{R}}\left\lbrace\exists \,P_M\,|\,P_M(x)=u(x)\;\;\mbox{and}\;\;P_M(y)\leq u(y),\;\;\forall y\in \mathcal{O}\right\rbrace.
\]
Finally, 
\[
	\Theta(u,\mathcal{O})(x)\,:=\,\max\left\lbrace \overline{\Theta}(u,\mathcal{O})(x),\,\underline{\Theta}(u,\mathcal{O})(x)\right\rbrace.
\]
\end{Definition}

\begin{Definition}\label{def_asegs}
Let $\mathcal{O}\subset\mathbb{R}^d$, $u\in\mathcal{C}(\mathcal{O})$ and $M>0$. We set
\[
	\overline{G}_M(u,\mathcal{O})\,:=\,\left\lbrace x\in \mathcal{O}\,|\,\exists \,P_M\;\mbox{s.t.}\;P_M(x)=u(x)\;\;\mbox{and}\;\;P_M(y)\geq u(y)\;\forall y\in \mathcal{O} \right\rbrace,
\]
\[
	\underline{G}_M(u,\mathcal{O})\,:=\,\left\lbrace x\in \mathcal{O}\,|\,\exists \,P_M\;\mbox{s.t.}\;P_M(x)=u(x)\;\;\mbox{and}\;\;P_M(y)\leq u(y)\;\forall y\in \mathcal{O} \right\rbrace
\]
and
\[
	G_M(u,\mathcal{O})\,:=\,\overline{G}_M(u,\mathcal{O})\,\cap\underline{G}_M(u,\mathcal{O}).
\]
In addition, we have
\[
\overline{A}_M(u,\mathcal{O})\,:=\,\mathcal{O}\setminus \overline{G}_M(u,\mathcal{O}),\;\;\;\;\;\;\;\;\;\;\underline{A}_M(u,\mathcal{O})\,:=\,\mathcal{O}\setminus \underline{G}_M(u,\mathcal{O})
\]
and
\[
A_M(u,\mathcal{O})\,:=\,\mathcal{O}\setminus G_M(u,\mathcal{O}).
\]
\end{Definition}

The next lemma relates the $L^p$-norms of the function $\Theta(u,\mathcal{O})$ with norms of $u$ in Sobolev spaces.

\begin{Lemma}\label{lem_prelimtheta1}
Let $p\in(d/2,+\infty)$. If $u\in W^{2,p}(\mathcal{O})$, we have
\[
	\left\|\Theta(u,\mathcal{O})\right\|_{L^p(\mathcal{O})}\,\leq\,C\left\|u\right\|_{W^{2,p}(\mathcal{O})},
\]
where $C=C(d,p)$ is a nonnegative constant.
\end{Lemma}
\begin{proof}
For the proof of Lemma \ref{lem_prelimtheta1}, we refer the reader to \cite[Lemma 2.5]{lizhang}.
\end{proof}

\begin{Lemma}\label{lem_Lpimpliesdecay}
Let $\mathcal{O}\subset\mathbb{R}^d$. Suppose $u\in W^{2,p}(\mathcal{O})$. Then, there exists $C>0$ such that
\[
	|A_t(u,\mathcal{O})|\,\leq\,Ct^{-p}.
\]
\end{Lemma}
\begin{proof}
Because $u\in W^{2,p}(\mathcal{O})$, Lemma \ref{lem_prelimtheta1} yields $\Theta(u,\mathcal{O})\in L^p(\mathcal{O})$. Therefore, there exists $C>0$ satisfying
\[
	|\left\lbrace x\in\mathcal{O}\,|\,\Theta(u,\mathcal{O})\,>\,t  \right\rbrace|\,\leq\, Ct^{-p}.
\]
Finally, notice that 
\[
	A_t(u,\mathcal{O})\,\subset\,\left\lbrace x\in\mathcal{O}\,|\,\Theta(u,\mathcal{O})\,>\,t  \right\rbrace.
\]
\end{proof}


We proceed by collecting elementary results on the Calder\'on-Zygmund decomposition and the maximal function; refer to \cite{ccbook}.

\begin{Lemma}[Calder\'on-Zygmund cube decomposition]\label{lem_czd}
Let $A\subset B\subset Q$ be measurable sets and $\sigma\in(0,1)$. For a dyadic cube $K$, denote by $\overline{K}$ its (unique) predecessor. If
\[
	|A|\,\leq\,\sigma
\]
and 
\[
	|A\cap K|\,\geq\,\sigma|K|\;\;\;\;\;\Rightarrow\;\;\;\;\;\overline{K}\subset B,
\]
then
\[
	|A|\,\leq\,\sigma|B|.
\]
\end{Lemma}

\begin{Lemma}\label{lem_maximal}
Let $h:\mathcal{O}\to\mathbb{R}$ be a measurable, nonnegative, function. For $t>0$, define
\[
	\mu_h(t)\,:=\,|\left\lbrace x\in\mathcal{O}\,|\, h(x)\,>\,t \right\rbrace|.
\]
Fix $\nu>0$ and $M>1$; for $p>1$, define
\[
	\mathcal{S}\,:=\,\sum_{k=1}^\infty M^{pk}\mu_h(\nu M^k).
\]
Then, $h\in L^p(\mathcal{O})$ if and only if $\mathcal{S}<\infty$. Moreover, there exists a positive constant $C=C(\nu,p,M)$ such that 
\[
	C^{-1}\mathcal{S}\,\leq\,\left\|h\right\|_{L^p(\mathcal{O})}^p\,\leq\,C\left(|\mathcal{O}|\,+\,\mathcal{S}\right).
\]
\end{Lemma}

In the sequel, we detail the main assumptions under which we work in this paper.

\subsection{Principal assumptions}\label{sec_assumptions}

In this section, we detail our main assumptions. We begin with the uniform ellipticity of the operators.

\begin{Assumption}[Uniform ellipticity]\label{assump_ellipticity}
We suppose the matrix $A_{\alpha,\beta}:B_1\times\mathcal{A}\times\mathcal{B}\to\mathbb{R}^{d^2}$ is $(\lambda,\Lambda)$-elliptic.
\end{Assumption}

Our next assumption describes the smallness regime required to prove Sobolev regularity.

\begin{Assumption}[Smallness regime -- Sobolev estimates]\label{assump_smallness}
We suppose that $\overline{A}_\beta:B_1\times\mathcal{B}\to\mathbb{R}^{d^2}$ satisfies
\[
	\left|A_{\alpha,\beta}(x)\,-\,\overline{A}_{\beta}(x)\right|\,\leq\,\epsilon_1
\]
uniformly in $x$, $\alpha$ and $\beta$, and
\[
	\left\|f\right\|_{L^p(B_1)}\,\leq\,\epsilon_1,
\]
where $\epsilon_1>0$ is a number to be determined further in the paper.
\end{Assumption}

Notice that A\ref{assump_ellipticity} and A\ref{assump_smallness} imply that $A_{\beta}$ is also $(\lambda,\Lambda)$-elliptic.
Now, we suppose that solutions to our limiting problem have estimates in $W^{2,q}$, for $q>d$ fixed.

\begin{Assumption}[$W^{2,q}$-estimates for the approximate problem]\label{assump_w2qest}
Let $r\in(0,1)$ and $d<p<q$. Let $h\in\mathcal{C}(B_r)$ be a viscosity solution to
\begin{equation}\label{eq_bellman}
	\inf_{\beta\in\mathcal{B}}\,\left[-\tr\left(\overline{A}_{\beta}(x)D^2h\right)\right]\,=\,0\;\;\;\;\;\mbox{in}\;\;\;\;\;B_r.
\end{equation}
Then, $h\in W^{2,q}(B_r)\cap\mathcal{C}(\overline{B}_{r+})$ and there exists a universal constant $C>0$ such that
\[
	\left\|h\right\|_{W^{2,q}(B_r)}\,\leq\,C.
\]
\end{Assumption}

It is reasonable to suppose that solutions to \eqref{eq_bellman} have estimates in $W^{2,q}(B_r)$. For example, this would follow from \cite[Theorem 7.1]{ccbook} provided the oscillation of $\overline{A}_\beta$ with respect to its fixed-coefficients counterpart is controlled in the $L^q$-sense. That is, the quantity 
\[
	\Theta_{x_0}(x)\,:=\,\sup_{\beta\in\mathcal{B}}\left|\overline{A}_\beta(x)\,-\,\overline{A}_\beta(x_0)\right|
\]
satisfies a smallness regime of the form
\[
	\left(\intav{B_r(x_0)}\left|\Theta_{x_0}(x)\right|^ddx\right)^\frac{1}{d}\,\ll\,1/2.
\]



We consider more general formulations of the Isaacs equation, when establishing estimates in Sobolev spaces. The following is an assumption on the coefficient of the lower-order terms.

\begin{Assumption}[Vector field $\bb_{\alpha,\beta}$]\label{assump_vf}
We suppose that $\bb_{\alpha,\beta}\in L^\infty(B_1)$ uniformly. That is, 
\[
	\sup_{\alpha\in\mathcal{A}}\sup_{\beta\in\mathcal{B}}\left\|\bb_{\alpha,\beta}\right\|_{L^\infty(B_1)}\,\leq\,C,
\]
for some $C>0$.
\end{Assumption}

We proceed by introducing further smallness regimes. As concerns the source term $f$, we make assumptions on its norms in BMO spaces. Let $x_0\in B_1$ and take $r_0:=\dist(x_0,\partial B_1)$. We denote by $\left\langle f \right\rangle_{x_0,r_0}$ the following quantity:
\[
	\left\langle f \right\rangle_{x_0,r_0}\,:=\,\intav{B_{r_0}(x_0)}\left|f(x)\right|dx;
\]
for convenience, we set $\left\langle f \right\rangle_{0,1}\equiv \left\langle f \right\rangle$.


The next assumption regards the smallness regime used in the study of $\mathcal{C}^{1,\llip}$-regularity for \eqref{eq_isaacs}.

\begin{Assumption}[Smallness regime -- estimates in $\mathcal{C}^{1,\llip}$]\label{assump_osc} 

{\color{white}.}

\begin{enumerate}

\item We suppose that, for every $x_0\in B_1$,
\[
	\sup_{x\in B_{r_0}}\,\left|A_{\alpha,\beta}(x)\,-\,\overline{A}_\beta(x_0)\right|\,\leq\,\epsilon_2,
\]
uniformly in $\alpha$ and $\beta$. 
\item Furthermore, the source term satisfies
\[
	\sup_{r\in(0,r_0]}\,\intav{B_{r}(x_0)}\left|f(x)\,-\,\left\langle f \right\rangle_{r_0,x_0}\right|^pdx\,\leq\,\epsilon_2^p,
\]
for every $x_0\in B_1$.
As before, $\epsilon_2>0$ will be determined further in the paper.
\end{enumerate}
\end{Assumption}

Lastly, we describe the smallness regime required to produce estimates in $\mathcal{C}^{2,\gamma}$ for the solutions of \eqref{eq_isaacs}.

\begin{Assumption}[Smallness regime -- estimates in $\mathcal{C}^{2,\gamma}$]\label{assump_oscc2alpha} We suppose that
\[
	\sup_{x\in B_r}\sup_{\alpha\in \mathcal{A}}\sup_{\beta\in \mathcal{B}}\left|A_{\alpha,\beta}(x)\,-\,\overline{A}_\beta\right|\,\leq\,\epsilon_3r^\gamma.
\]
Moreover, the source term satisfies
\[
	\intav{B_r}|f(x)|^pdx\,\leq\,\epsilon_3^pr^{\gamma p}.
\]
As before, $\epsilon_3>0$ will be determined further in the paper.
\end{Assumption}


\section{Estimates in Sobolev spaces}\label{isaacsw2p}

In this section, we present the proof of Theorem \ref{teo_main1}. We start with a result on the regularity in $W^{2,p}$ for the solutions to \eqref{eq_isaacs}. That is, the equation without dependence on the gradient. In this setting, we prove the following proposition.

\begin{Proposition}\label{prop_w2peasy}
Let $u\in\mathcal{C}(B_1)$ be a viscosity solution to \eqref{eq_isaacs}. Let $d<p<q$. Suppose A\ref{assump_ellipticity}-A\ref{assump_w2qest} are in force. Then, $u\in W^{2,p}_{loc}(B_1)$ and
\[
	\left\|u\right\|_{W^{2,p}(B_{1/2})}\,\leq\,C\left(\left\|u\right\|_{L^\infty(B_1)}\,+\,\left\|f\right\|_{L^p(B_1)}\right),
\]
where $C>0$ is a universal constant.
\end{Proposition}

After that, we argue as in \cite{swiech97} and prove $W^{1,\overline{p}}$  estimates for \eqref{eq_isaacslot}. 

\begin{teo}[Estimates in $W^{1,\overline{p}}$]\label{teo_w1p}
Let $u\in\mathcal{C}(B_1)$ be a viscosity solution to \eqref{eq_isaacslot}. Suppose A\ref{assump_ellipticity}-A\ref{assump_vf} are in force. Then, for every $\overline{p}\in(1,\infty)$, we have $u\in W^{1,\overline{p}}_{loc}(B_1)$ and there exists $C>0$ such that
\[
	\left\|u\right\|_{W^{1,\overline{p}}(B_{1/2})}\,\leq\,C\left(\left\|u\right\|_{L^{\infty}(B_{1})}\,+\,\left\|f\right\|_{L^{p}(B_{1})}\right).
\]
\end{teo}

To conclude the proof of Theorem \ref{teo_main1}, we resort to a reduction argument. Next, we detail the proof of Proposition \ref{prop_w2peasy}.

\subsection{Proof of Proposition \ref{prop_w2peasy}}\label{sec_proofw2peasy}

We proceed with an estimate in $W^{2,\mu}$.
\begin{Lemma}\label{lem_w2delta}
Let $u\in\mathcal{C}(B_1)$ be a viscosity solution to \eqref{eq_isaacs}. Then, there exist universal constants $C>0$ and $\mu>0$ such that
\begin{equation}\label{eq_decayw2delta}
	|A_t(u,B_1)\,\cap\,Q|\,\leq\,Ct^{-\mu},
\end{equation}
for every $t>0$.
\end{Lemma}

The first version of Lemma \ref{lem_w2delta} has appeared in \cite{fanghuaw2delta} and addressed the linear case. In \cite[Lemma 7.4]{ccbook} the authors have proven the result for fully nonlinear elliptic operators. Next, an approximation lemma relates the solutions of \eqref{eq_isaacs} to an auxiliary function. Regularity properties for the latter imply a refinement of the decay rate in \eqref{eq_decayw2delta}.

\begin{Proposition}[First Approximation Lemma]\label{lem_prox1}
Let $u\in\mathcal{C}(B_1)$ be a viscosity solution to \eqref{eq_isaacs}. Suppose A\ref{assump_ellipticity}-A\ref{assump_w2qest} are in force. Then, there exists $h\in W^{2,q}(B_{7/8})\cap\mathcal{C}(\overline{B}_{8/9})$ such that
\[
	\left\|h\right\|_{W^{2,q}(B_{7/8})}\,\leq\,C
\]
and
\[
	\left\|u\,-\,h\right\|_{L^\infty(B_{7/8})}\,\leq\,C_1\left(\epsilon_1^\gamma\,+\,\left\|f\right\|_{L^\infty(B_1)}\right),
\]
where $C_1=C_1(d,\lambda,\Lambda,C,q,p)$ and $\gamma=\gamma(d,\lambda,\Lambda,C,q,p)$ are nonnegative constants. 
\end{Proposition}
\begin{proof}
Consider the function $h\in\mathcal{C}(\overline{B}_{8/9})$, solution to
\begin{equation*}
	\begin{cases}
		\inf_{\beta\in\mathcal{B}}\left[-\tr\left(\overline{A}_{\beta}(x)D^2h\right)\right]\,=\,0&\;\;\;\;\;\mbox{in}\;\;\;\;\;B_{8/9}\\
		h\,=\,u&\;\;\;\;\;\mbox{on}\;\;\;\;\;\partial B_{8/9}.
	\end{cases}
\end{equation*} 
Because of A\ref{assump_w2qest}, we have
\[
	\left\|h\right\|_{W^{2,q}(B_{7/8})}\,\leq\,C.
\]
In addition, standard results in interior H\"older regularity imply
\begin{equation}\label{eq_ualpha}
	\left\|u\right\|_{\mathcal{C}^{\overline{\gamma}}(\overline{B}_{8/9})}\,\leq\,C\left(\left\|u\right\|_{L^\infty(B_1)}\,+\,\left\|f\right\|_{L^p(B_1)}\right)
\end{equation}
and
\begin{equation}\label{eq_halpha}
	\left\|h\right\|_{\mathcal{C}^\frac{\overline{\gamma}}{2}(\overline{B}_{8/9})}\,\leq\,C\left(\left\|u\right\|_{L^\infty(B_1)}\,+\,\left\|f\right\|_{L^p(B_1)}\right),
\end{equation}
for some $\overline{\gamma}\in(0,1)$.

Next, we fix $\delta\in(0,1/2)$ and take $x_0\in B_{8/9-\delta}$; clearly $B_\delta(x_0)\subset B_{8/9}$. Take $x_1\in\partial B_\delta(x_0)$; by applying $W^{2,q}$-estimates to the function $h(x)-h(x_1)$ in $B_\delta(X_0)$, we conclude
\begin{align*}\label{eq_hessianx0}
	\left\|D^2h\right\|_{L^q(B_{\delta/2}(x_0))}\,&\leq\,C_h\delta^{\frac{d-2q}{q}}\left\|h\,-\,h(x_1)\right\|_{L^\infty(\partial B_{\delta}(x_0))}\\
		&\leq\,C\delta^{\frac{d-2q}{q}+\frac{\overline{\gamma}}{2}}\left(\left\|u\right\|_{L^\infty(B_1)}\,+\,\left\|f\right\|_{L^p(B_1)}\right),
\end{align*}
where the second inequality follows from \eqref{eq_halpha}; therefore,
\begin{equation*}\label{eq_hessianlp}
	\left\|D^2h\right\|_{L^q(B_{8/9-\delta})}\,\leq\,C\delta^{\frac{d-2q}{q}+\frac{\overline{\gamma}}{2}-d}\left(\left\|u\right\|_{L^\infty(B_1)}\,+\,\left\|f\right\|_{L^p(B_1)}\right).
\end{equation*}
By A\ref{assump_smallness} we obtain
\[
	\left|\sup_{\alpha\in\mathcal{A}}\inf_{\beta\in\mathcal{B}}\left[-\tr\left(A_{\alpha,\beta}(x_0)D^2h(x_0)\right)\right]\right|\leq C(d)\epsilon_1\left\|D^2h(x_0)\right\|,
\]
for $x_0\in B_{8/9-\delta}$. Then, 
\begin{align}\label{eq_hessiangeral}
	&\nonumber\left\|\sup_{\alpha\in\mathcal{A}}\inf_{\beta\in\mathcal{B}}\left[-\tr\left(A_{\alpha,\beta}(x)D^2h(x)\right)\right]\right\|_{L^p(B_{8/9-\delta})}\leq C(d)\epsilon_1\left\|D^2h\right\|_{L^q(B_{8/9-\delta})}\\
	&\qquad\leq C(d)\epsilon_1\delta^{\frac{d-2q}{q}+\frac{\overline{\gamma}}{2}-d}\left(\left\|u\right\|_{L^\infty(B_1)}\,+\,\left\|f\right\|_{L^p(B_1)}\right).
\end{align}

We combine \eqref{eq_ualpha} and \eqref{eq_halpha} to obtain
\begin{equation}\label{eq_umenosh}
	\left\|u\,-\,h\right\|_{L^\infty(B_{8/9-\delta})}\,\leq\,C\delta^\frac{\overline{\gamma}}{2}\left(\left\|u\right\|_{L^\infty(B_1)}\,+\,\left\|f\right\|_{L^p(B_1)}\right).
\end{equation}

By gathering the information in \eqref{eq_hessiangeral} and \eqref{eq_umenosh} and using the maximum principle, we get
\begin{align*}
	\left\|u\,-\,h\right\|_{L^\infty(B_{8/9-\delta})}&\leq C\delta^\frac{\overline{\gamma}}{2}\left(\left\|u\right\|_{L^\infty(B_1)}\,+\,\left\|f\right\|_{L^p(B_1)}\right)+C\left\|f\right\|_{L^p(B_{8/9-\delta})}\\
	&\quad +\, C(d)\epsilon_1\delta^{\frac{d-2q}{q}+\frac{\overline{\gamma}}{2}-d}\left(\left\|u\right\|_{L^\infty(B_1)}\,+\,\left\|f\right\|_{L^p(B_1)}\right)\\
	&\leq \left(C\delta^\frac{\overline{\gamma}}{2}+C(d)\epsilon_1\delta^{\frac{d-2q}{q}+\frac{\overline{\gamma}}{2}-d}\right)\left(\left\|u\right\|_{L^\infty(B_1)}\,+\,\left\|f\right\|_{L^p(B_1)}\right)\\
	&\quad+\,C\left\|f\right\|_{L^p(B_{1})}.
\end{align*}
To conclude the proof, we set 
\[
	\delta\,:=\,\epsilon_1^{\frac{q}{dq\,+\,2q\,-\,d}}\;\;\;\;\;\;\;\;\;\;\;\;\mbox{and}\;\;\;\;\;\;\;\;\;\;\;\;\gamma\,:=\,\frac{\overline{\gamma} q}{2(dq\,+\,2q\,-\,d)}
\]
and observe that $B_{7/8}\subset B_{8/9-\delta}$. 
\end{proof}


\begin{Lemma}\label{lem_refine1}
Let $u\in\mathcal{C}(B_1)$ be a viscosity solution to \eqref{eq_isaacs}. Suppose that A\ref{assump_ellipticity}-A\ref{assump_w2qest} are in force. Suppose further that
\[
	-|x|^2\,\leq\,u(x)\,\leq\,|x|^2\;\;\;\;\;\mbox{in}\;\;\;\;\;B_1\setminus B_{6/7}.
\]
Then, there exist $\overline{M}>1$ and $\sigma\in(0,1)$ such that
\[
	|G_{\overline{M}}(u,B_1)\cap Q|\,\geq 1\,-\,\sigma.
\]
\end{Lemma}
\begin{proof}
Let $h$ be the approximate function from Proposition \ref{lem_prox1}. Extend $h$ outside $B_{7/8}$ to have
\[
	h\,\equiv\, u\;\;\;\;\;\mbox{in}\;\;\;\;\;B_1\setminus B_{8/9}
\]
and
\[
	\left\|u\,-\,h\right\|_{L^\infty(B_1)}\,=\,\left\|u\,-\,h\right\|_{L^\infty(B_{7/8})}.
\]
A standard application of the maximum principle yields $\left\|h\right\|_{L^\infty(B_{7/8})}\leq 1$. Therefore, 
\begin{equation}\label{eq_hw2qfora}
	-2\,-\,|x|^2\,\leq\,h(x)\,\leq\,|x|^2\,+\,2\;\;\;\;\;\mbox{in}\;\;\;\;\;B_1\setminus B_{7/8}.
\end{equation}
In addition, we conclude from Proposition \ref{lem_prox1} that $h\in W^{2,q}(B_{7/8})$. By combining this fact with \eqref{eq_hw2qfora}, we secure the existence of a constant $C>0$ for which
\begin{equation}\label{eq_Ahw2q}
	|A_t(h,B_1)\cap Q|\,\leq\,Ct^{-q}.
\end{equation}

Next, define the auxiliary function 
\[
	v\,:=\,\frac{u\,-\,h}{2C_1\epsilon_1^\gamma}.
\]
We have that $$v\in S\left(\lambda,\Lambda,\frac{f}{2C_1\epsilon_1^\gamma}\right)$$ satisfies the assumptions of Lemma \ref{lem_w2delta}. Hence,
\[
	|A_t(u-h,B_1)\cap Q|\,\leq\,C\epsilon_1^{\gamma\mu}t^{-\mu}.
\]
Write $\overline{M}\equiv 2t$ to obtain
\begin{align*}
	|A_{\overline{M}}(u,B_1)\cap Q|\,&\leq\, |A_{\overline{M}/2}(u-h,B_1)\cap Q|\,+\,|A_{\overline{M}/2}(h,B_1)\cap Q|\\
	&\leq C\epsilon_1^{\gamma\mu}\left(\frac{\overline{M}}{2}\right)^{-\mu}\,+\,\overline{C}\left(\frac{\overline{M}}{2}\right)^{-q}.
\end{align*}
By taking $\epsilon_1$ small enough and choosing $\overline{M}$ appropriately, we conclude the proof.
\end{proof}

We notice that, because $d<p<q$, there will always exist a constant $\overline{M}$ allowing us to conclude the proof of Lemma \ref{lem_refine1}. This is a rigorous instance where the (intuitive) requirement $p<q$ is binding for the theory.

\begin{Lemma}\label{lem_refine2}
Let $u\in\mathcal{C}(B_1)$ be a viscosity solution to \eqref{eq_isaacs}. Suppose A\ref{assump_ellipticity}-A\ref{assump_w2qest} hold. Let $\overline{Q}$ be a cube such that $Q\subset \overline{Q}$. If
\[
	G_1(u,B_1)\cap \overline{Q}\,\neq\,\emptyset,
\]
we have
\[
	|G_{M}(u,B_1)\cap Q|\,\geq 1-\sigma,
\]
for some $M>1$.
\end{Lemma}
\begin{proof}
Because $G_1(u,B_1)\cap \overline{Q}$ is nonempty, there exists $x_0\in B_1$ and an affine function $\ell(x)$ such that
\[
	-\frac{|x\,-\,x_0|^2}{2}\,\leq\,u(x)\,-\,\ell(x)\,\leq\,\frac{|x\,-\,x_0|^2}{2}
\]
in $B_1$. Next, choose $C_3>0$ to guarantee that
\[
	v(x)\,:=\,\frac{(u\,-\,\ell)(x)}{C_3}
\]
satisfies $\left\|v\right\|_{L^\infty(B_1)}\leq 1$ and
\[
	-|x|^2\,\leq\,v(x)\,\leq\,|x|^2\;\;\;\;\;\;\;\;\;\;\mbox{in}\;\;\;\;\;\;\;\;\;\;B_1\setminus B_{6/7}.
\]
Therefore, Lemma \ref{lem_refine1} yields
\[
	|G_{\overline{M}}(v,B_1)\cap Q|\,\geq\,1\,-\,\sigma.
\]
Set $M:= C_3\,\overline{M}$; we conclude that 
\[
	|G_M(u,B_1)\cap Q|\,=\,|G_{C_3\overline{M}}(u,B_1)\cap Q|\,=\,|G_{\overline{M}}(v,B_1)\cap Q|\,\geq\,1\,-\,\sigma,
\]
which finishes the proof.
\end{proof}
We note that, in the remainder of this section, the constant $M$ refers to $M:=C_3\overline{M}$.

\begin{Lemma}\label{lem_refine3}
Let $u\in\mathcal{C}(B_1)$ be a viscosity solution to \eqref{eq_isaacs}. Suppose A\ref{assump_ellipticity}-A\ref{assump_w2qest} are in force. Extend $f$ by zero outside of $B_1$ and set
\[
	A\,:=\,A_{M^{k+1}}(u,B_1)\cap Q
\]
and
\[
	B\,:=\,\left(A_{M^k}(u,B_1)\cap Q\right)\,\cup\,\left\lbrace x\in Q\,|\,m\left(f^p\right)\,\geq\,\left(CM^k\right)^p \right\rbrace.
\]
Then, 
\[
	|A|\,\leq\,\sigma|B|.
\]
\end{Lemma}
\begin{proof}
In what follows, we make use of the Calder\'on-Zygmund decomposition. Lemma \ref{lem_refine1} yields
\[
	|G_{M^{k+1}}(u,B_{1})\cap Q|\geq |G_{M}(u,B_{1})\cap Q|\geq |G_{C_3\overline{M}}(u,B_{1})\cap Q|\geq 1-\sigma. 
\]
The definition of $A$ implies then
\[
	|A|\,\leq\,\sigma.
\]

Set $K:=Q_{1/2^i}(x_0)$ and denote by $\overline{K}$ the predecessor of $K$. To conclude the proof we must verify
\begin{equation}\label{eq_cont3}
	\left(|A_{M^{k+1}}(u,B_{1})\cap K|\,=\,|A\cap K|\,>\,\sigma |K|\right)\;\;\;\;\;\;\Rightarrow\;\;\;\;\;\;\;\overline{K}\subset B.
\end{equation}
We proceed by contradiction; suppose $\overline{K}\not\subset B$. Then, there exists $x_1$ satisfying both
\begin{equation}\label{eq_cont1}
	x_1\,\in\,\overline{K}\cap G_{M^k}(u,B_{1})
\end{equation}
and
\begin{equation}\label{eq_cont2}
	m(f^p)(x_1)\,<\,(CM^k)^p.
\end{equation}
We introduce an affine transformation $T:Q\to K$, given as follows:
\[
	T(y)\,:=\,x_0\,+\,\frac{y}{2^i}.
\]
Then, consider the auxiliary function
\[
	\overline{u}(y)\,:=\,\frac{2^{2i}}{M^k}(u\circ T)(y)\,=\,\frac{2^{2i}}{M^k}u\left(x_0\,+\,\frac{y}{2^i}\right).
\]
We proceed by verifying that $\overline{u}$ satisfies the assumptions of Lemma \ref{lem_refine2}. Notice that
\[
	D^2\overline{u}(y)\,=\,\frac{1}{M^k}D^2u\left(x_0\,+\,\frac{y}{2^i}\right);
\]
therefore, 
\[
	\sup_{\alpha\in\mathcal{A}}\inf_{\beta\in\mathcal{B}}\,\left[-\tr\left(A_{\alpha,\beta}\left(x+\frac{y}{2^i}\right)D^2\overline{u}\right)\right]\,=\,\overline{f}\;\;\;\;\;\mbox{in}\;\;\;\;\;B_1,
\]
where
\[
	\overline{f}(y)\,:=\,\frac{f(x_0+2^{-i}y)}{M^k}.
\]
We have
\[
	\left\|\overline{f}\right\|^p_{L^p(B_{2^{-i}})}\,=\,\frac{2^{i(2-d)}}{M^{kp}}\int_{B_{2^{-i}}(x_0)}|f(x)|^pdx\,\leq\,\epsilon_1^p,
\]
where the last inequality stems from \eqref{eq_cont2}. From \eqref{eq_cont1} we infer that
\[
	G_1(\overline{u},B_{2^{-i}}(x_0))\cap\overline{Q}\neq\emptyset.
\]
As a consequence, we obtain
\[
	|G_M(\overline{u},B_{2^{-i}}(x_0))\cap Q|\,\geq \,(1-\sigma)|Q|.
\]
That is,
\[
	|G_{M^{k+1}}(u,B_1)\cap K|\,\geq\,(1-\sigma)|K|,
\]
which contradicts \eqref{eq_cont3} and finishes the proof.
\end{proof}

Next, we complete the proof of Proposition \ref{prop_w2peasy}.

\begin{proof}[Proof of Proposition \ref{prop_w2peasy}]
To prove the proposition, it suffices to verify that
\begin{equation}\label{eq_summation}
	\sum_{k=1}^\infty M^{pk}|A_{M^k}(u,B_{1/2})|\,\leq\,C,
\end{equation}
for some constant $C>0$. Consider the quantities
\[
	a_k\,:=\,|A_{M^k}(u,B_1)\cap Q|\;\;\;\;\;\mbox{and}\;\;\;\;\;b_k\,:=\,\left|\left\lbrace x\in Q\,|\,m(f^p)(x)\,\geq\,(CM^k)^p \right\rbrace\right|.
\]
The following inequality is due to Lemma \ref{lem_refine3}:
\[
	a_{k+1}\,\leq\,a_k\,+\,b_k.
\]
Hence,
\[
	a_k\,\leq\,\sigma^k\,+\,\sum_{i=0}^{k-1}\sigma^{k-i}b_i.
\]
It is clear that $f^q\in L^{p/q}(B_1)$; as regards the maximal function, it implies $m(f^q)\in L^{p/q}(B_1)$ and
\[
	\left\|m(f^q)\right\|_{L^{p/q}(B_1)}\,\leq\,C\left\|f\right\|_{L^p(B_1)}^q\,\leq\,C,
\]
for some $C>0$. Elementary properties of the maximal functions imply
\[
	\sum_{k=0}^\infty M^{pk}b_k\,\leq\,C.
\]
We have found that
\begin{align*}
	\sum_{k=1}^\infty M^{pk}a_k\,&\leq\,\sum_{k=1}^\infty(\sigma M^p)^k\,+\,\sum_{k=1}^\infty\sum_{i=0}^{k-1}\sigma^{k-i}M^{p(k-i)}M^{pi}b_i\\
		&\leq\,\sum_{k=1}^\infty(\sigma C_3\overline{M}^p)^k\,+\,\sum_{k=1}^\infty\sum_{i=0}^{k-1}\sigma^{k-i}M^{p(k-i)}M^{pi}b_i\\
		&\leq\,\sum_{k=1}^\infty 2^{-k}\,+\,\left(\sum_{i=0}^\infty M^{pi}b_i\right)\left(\sum_{j=1}^\infty 2^{-j}\right)\\
		&\leq\,C.
\end{align*}
This finishes the proof.
\end{proof}

Now, we continue with the proof of Theorem \ref{teo_w1p}. We follow closely the arguments in \cite{swiech97}.

\subsection{Proof of Theorem \ref{teo_w1p}}\label{sec_proofw1p}

We start with an approximation lemma:

\begin{Proposition}[Second Approximation Lemma]\label{prop_approxw1p}
Let $u\in\mathcal{C}(B_1)$ be a viscosity solution to \eqref{eq_isaacslot}. Suppose A\ref{assump_ellipticity}-A\ref{assump_w2qest} and A\ref{assump_vf} are in force. For every $\delta>0$ it is possible to choose $\varepsilon_1>0$ so that there exists $h\in W^{2,q}_{loc}(B_1)$ satisfying
\begin{equation}
		\begin{cases}
			\inf_{\beta\in\mathcal{B}}\left[-\tr\left(A_\beta(x)D^2h(x)\right)\right]\,=\,0&\;\;\;\;\;\mbox{in}\;\;\;\;\;B_{8/9}\\
			h\,=\,u&\;\;\;\;\;\mbox{on}\;\;\;\;\;\partial B_{8/9},
		\end{cases}
\end{equation}
with 
\begin{equation}
	\left\|h\right\|_{W^{2,q}(B_{1/2})}\,\leq\,C
\end{equation}
and
\[
	\left\|u\,-\,h\right\|_{L^\infty(B_{7/8})}\,\leq\,\delta,
\]
for some $C>0$, universal.
\end{Proposition}
\begin{proof}
We prove the proposition by contradiction. Suppose its statement is false. There would be a sequence $(A_{\alpha,\beta}^n)_{n\in\mathbb{N}}$ of matrices $A_{\alpha,\beta}^n:B_1\times\mathcal{A}\times\mathcal{B}\to\mathbb{R}^{d^2}$, a sequence $(\bb_{\alpha,\beta}^n)_{n\in\mathbb{N}}$ of vector fields $\bb_{\alpha,\beta}^n:B_1\times\mathcal{A}\times\mathcal{B}\to\mathbb{R}^{d}$ and sequences of functions $(u_n)_{n\in\mathbb{N}}$ and $(f_n)_{n\in\mathbb{N}}$ satisfying
\begin{equation}\label{eq_auxw2phardcore}
	\sup_{\alpha\in\mathcal{A}}\inf_{\beta\in\mathcal{B}}\,\left[-\tr\left(A_{\alpha,\beta}^n(x)D^2u_n\right)\right]\,=\,f_n\;\;\;\;\;\mbox{in}\;\;\;\;\;B_1,
\end{equation}
with
\[
	\left|A_{\alpha,\beta}^n(x)\,-\,\overline{A}_\beta(x)\right|\,+\,\left\|\bb_{\alpha,\beta}^n\right\|_{L^\infty(B_1)}\,+\,\left\|f_n\right\|_{L^p(B_1)}\,<\,\frac{1}{n}
\]
such that
\[
	\left\|u_n\,-\,h\right\|_{L^\infty(B_{7/8})}\,>\,\delta_0
\]
for some $\delta_0>0$ and every solution $h$ to
\[
	\inf_{\beta\in\mathcal{B}}\,\left[-\tr\left(\overline{A}_{\beta}(x)D^2h\right)\right]\,=\,0\;\;\;\;\;\mbox{in}\;\;\;\;\;B_{8/9}.
\]
The regularity theory available for \eqref{eq_auxw2phardcore} implies that $u_n$ converges, through a subsequence if necessary, to a function $u_\infty$ in the $\mathcal{C}^\nu$-topology; see, \cite[Lemma 1.9]{swiech97}. Standard results on the stability of viscosity solutions yield
\[
	\inf_{\beta\in\mathcal{B}}\,\left[-\tr\left(\overline{A}_{\beta}(x)D^2u_\infty\right)\right]\,=\,0\;\;\;\;\;\mbox{in}\;\;\;\;\;B_{8/9};
\]
see, for example, \cite[Lemma 1.7]{swiech97}. Because of A\ref{assump_w2qest}, we have $u_\infty\in W^{2,q}_{loc}(B_1)$. By setting $h\equiv u_\infty$ we obtain a contradiction and complete the proof.
\end{proof}

In the sequel, Proposition \ref{prop_approxw1p} builds upon scaling properties of \eqref{eq_isaacslot} to produce integral estimates for the gradient of solutions.

\begin{proof}[Proof of Theorem \ref{teo_w1p}]
The proof is similar to the one in \cite[Theorem 2.1]{swiech97}. We omit the details here.

\end{proof}

Theorem \ref{teo_w1p} is instrumental in proving Theorem \ref{teo_main1}. It shows that estimates in $W^{1,p}$ can be accessed through approximations by a Belllman operator.

\subsection{Proof of Theorem \ref{teo_main1}}\label{sec_proofteo_main1}

Here, we conclude the proof of Theorem \ref{teo_main1}.

\begin{proof}[Proof of Theorem \ref{teo_main1}]
We know that $u$ solves \eqref{eq_isaacslot} pointwise a.e. in $B_1$; see \cite{CCKS}. Set
\[
	g(x)\,:=\,\sup_{\alpha\in\mathcal{A}}\inf_{\beta\in\mathcal{B}}\left[-\tr\left(A_{\alpha,\beta}(x)D^2u(x)\right)\right].
\]
We have
\[
|g(x)|\,\leq\, \sup_{\alpha\in\mathcal{A}}\sup_{\beta\in\mathcal{B}}\left|\bb_{\alpha,\beta}(x)\right||Du(x)|\,+\,|f(x)|\,\in\,L^p(B_1),
\]
because of Theorem \ref{teo_w1p}. Standard results on the equivalence of solutions imply that $u$ solves
\[
	\sup_{\alpha\in\mathcal{A}}\inf_{\beta\in\mathcal{B}}\left[-\tr\left(A_{\alpha,\beta}(x)D^2u(x)\right)\right]\,=\,g(x)\;\;\;\;\;\mbox{in}\;\;\;\;\;B_1
\]
in the viscosity sense. We refer the reader to \cite[Corollary 1.6]{swiech97}; see also \cite[Theorem 3.3]{ckss96}.
The result follows from Proposition \ref{prop_w2peasy}.
\end{proof}

\begin{Remark}\label{rem_bmo}
Once regularity in $W^{2,p}$ is available, it is possible to produce estimates for the Hessian of the solutions in John-Nirenberg spaces. It amounts to establishing the existence of a universal constant $C>0$ for which 
\[
	\sup_{r}\intav{B_r}\left|D^2u\,-\,\left\langle D^2u\right\rangle_r\right|^pdx\,\leq\,C\left(\left\|u\right|_{L^\infty(B_1)}\,+\,\left\|f\right\|_{L^p(B_1)}\right).
\]
The former estimate follows from two main ingredients. First, notice that Proposition \ref{prop_approxw1p} can be easily adapted to produce a value of $\varepsilon_1>0$ and a constant $0<\nu\ll 1/2$ so that
\[
	\sup_{x\in B_\nu}|u(x)\,-\,P_\nu(x)|\,<\,\nu^2,
\]
where $P_\nu$ is a paraboloid satisfying $\|P_\nu\|\leq C$. Then, an induction argument builds upon the previous inequality to prove the existence of a sequence of approximating polynomials. We refer the reader to \cite{silvtei}, \cite{teixeirauniversal} and \cite{pimtei} for further details.
\end{Remark}

\begin{Remark}\label{rem_c1alpha}
We observe that Theorem \ref{teo_main1} also implies estimates for the solutions to \eqref{eq_isaacslot} in H\"older spaces. Since $D^2u\in L^p_{loc}(B_1)$ for $d<p<q$, it yields $u\in\mathcal{C}^{1,\gamma^*}_{loc}(B_1)$, with
\[
	\gamma^*\,:=\,\min\left\lbrace \gamma_0,\,1\,-\,\frac{d}{q}\right\rbrace,
\]
where $\gamma_0$ is the exponent from the Krylov-Safonov theory. Therefore, the regularity of the approximate problem prevents $\gamma^*=\min\left\lbrace \gamma_0,1^-\right\rbrace$. \end{Remark}

\section{Estimates in $\mathcal{C}^{1,\llip}_{loc}(B_1)$}\label{isaacsc1ll}

Next, we examine the regularity of solutions to \eqref{eq_isaacs} in $\mathcal{C}^{1,\llip}_{loc}(B_1)$. We start with an approximation lemma. 

\begin{Proposition}\label{prop_prox2}
Let $u\in\mathcal{C}(B_1)$ be a viscosity solution to \eqref{eq_isaacs}. Suppose A\ref{assump_ellipticity} and A\ref{assump_osc}.1 are in force and $\left\|f\right\|_{L^p(B_1)}\ll 1$. 
For every $\delta>0$, it is possible to choose $\epsilon_2=\epsilon_2(\delta)>0$ to ensure the existence of $h\in\mathcal{C}^{2,\gamma}_{loc}(B_{1})$ satisfying 
\begin{equation}\label{eq_bellmanc1ll}
	\begin{cases}
		\inf_{\beta\in\mathcal{B}}\,\left[-\tr\left(\overline{A}_\beta(0)D^2h\right)\right]\,=\,0&\;\;\;\;\;\mbox{in}\;\;\;\;\;B_{8/9}\\
		h\,=\,u&\;\;\;\;\;\mbox{in}\;\;\;\;\;\partial B_{8/9},
	\end{cases}
\end{equation}
with
\begin{equation}\label{eq_hholder}
	\left\|h\right\|_{\mathcal{C}^{2,\gamma}(B_{1/2})}\,\leq\,C
\end{equation}
and
\[
	\left\|u\,-\,h\right\|_{L^\infty(B_{7/8})}\,\leq\,\delta.
\]
Moreover, the constants $C>0$ and $\gamma\in(0,1)$ are universal.
\end{Proposition}
\begin{proof}
Except for minor modifications, the result follows along the same lines as in the proof of Proposition \ref{prop_approxw1p}.
\end{proof}

\begin{Proposition}\label{prop_c1ll1}
Let $u\in\mathcal{C}(B_1)$ be a viscosity solution to \eqref{eq_isaacs}. Suppose A\ref{assump_ellipticity} and A\ref{assump_osc} are in force. 
There exists $0< \rho\ll 1$ and a sequence of polynomials $(P_n)_{n\in\mathbb{N}}$, given by
\[
	P_n(x)\,:=\,a_n\,+\,\bb_n\cdot x\,+\,\frac{1}{2}x^TC_nx,
\]
satisfying:
\begin{equation}\label{eq_condc1ll1}
	\inf_{\beta\in\mathcal{B}}\,\tr\left(\overline{A}_\beta(0) \,C_n\right)\,=\,\left\langle f\right\rangle,
\end{equation}
\begin{equation}\label{eq_condc1ll2}
	\sup_{B_{\rho^n}}\,|u(x)\,-\,P_n(x)|\,\leq\,\rho^{2k}
\end{equation}
and
\begin{equation}\label{eq_condc1ll3}
	|a_n\,-\,a_{n-1}|\,+\,\rho^{n-1}|\bb_n\,-\,\bb_{n-1}|\,+\,\rho^{2(n-1)}|C_n\,-\,C_{n-1}|\,\leq\, C\rho^{2(n-1)},
\end{equation}
for all $n\geq 0$.
\end{Proposition}
\begin{proof}
For ease of presentation, we split the proof in four steps.

\bigskip

\noindent{\bf Step 1}

\bigskip

We prove the result by induction in $n\in\mathbb{N}$. Set $P_{-1}=P_0=(1/2)x^TQx$, where $Q$ is such that
\[
	\inf_{\beta\in\mathcal{B}}\,\left[-\tr\left(\overline{A}_\beta(0)Q\right)\right]\,=\,\left\langle f\right\rangle.
\]
The case $n=0$ is clear. Suppose the case $n=k$ has been verified. We consider the case $n=k+1$. Define an auxiliary function $v_k:B_1\to\mathbb{R}$ as
\[
	v_k(x)\,:=\,\frac{(u\,-\,P_k)(\rho^kx)}{\rho^{2k}}.
\]
Notice that 
\[
	D^2v_k(x)\,=\,D^2u(\rho^kx)\,-\,C_k.
\]
Therefore, $v_k$ solves
\begin{equation}\label{eq_auxc1ll2}
	\sup_{\alpha\in\mathcal{A}}\inf_{\beta\in\mathcal{B}}\,\left[-\tr\left(A_{\alpha,\beta}(\rho^kx)(D^2v_k\,+\,C_k)\right)\right]\,=\,f_k\;\;\;\;\;\mbox{in}\;\;\;\;\;B_1,
\end{equation}
and 
\[
	f_k(x)\,:=\, f(\rho^kx).
\]

\bigskip

\noindent{\bf Step 2}

\bigskip

In what follows, we examine \eqref{eq_auxc1ll2}. Our goal is to approximate $v_k$ by suitable functions. We use A\ref{assump_osc} to conclude
\begin{align*}
	\left\|f_k\,-\,\left\langle f_k\right\rangle\right\|_{L^p(B_1)}^p\,&=\,\frac{1}{|\rho^{kd}|}\int_{B_{\rho^k}}\left|f(y)\,-\,\left\langle f\right\rangle_{\rho^k}\right|^pdy\\&\leq\,\sup_{r\in(0,1]}\,\intav{B_{r}}\left|f(y)\,-\,\left\langle f \right\rangle_r\right|^pdy\\&\leq\,\epsilon_2^p.
\end{align*}
Furthermore, we note that
\[
	\left|A_{\alpha,\beta}(\rho^kx)\,-\,\overline{A}_{\beta}(0)\right|\,\leq\,\epsilon_2.
\]
Consider next a function $v\in \mathcal{C}(B_1)$, viscosity solutions to 
\[
	\inf_{\beta\in\mathcal{B}}\,\left[-\tr\left(\overline{A}_\beta(0)D^2v\right)\right]\,=\,0;
\]
due to the Evans-Krylov theory, we have $v\in\mathcal{C}^{2,\gamma}_{loc}(B_1)$, with 
\[
	\left\|v\right\|_{\mathcal{C}^{2,\gamma}(B_{1/2})}\,\leq C^*.
\]
However, 
\[
	\inf_{\beta\in\mathcal{B}}\,\left[-\tr\left(\overline{A}_{\beta}(0)C_k\right)\right]\,=\,\left\langle f \right\rangle.
\]
Therefore, solutions to
\[
	\inf_{\beta\in\mathcal{B}}\,\left[-\tr\left(\overline{A}_{\beta}(0)(D^2h\,+\,C_k)\right)\right]\,=\,\left\langle f\right\rangle \;\;\;\;\;\mbox{in}\;\;\;\;\;B_1	
\]
also satisfy $h\in\mathcal{C}^{2,\gamma}_{loc}(B_1)$, with 
\[
	\left\|h\right\|_{\mathcal{C}^{2,\gamma}(B_{1/2})}\,\leq C=C(\left\langle f\right\rangle, C^*).
\]

\bigskip

\noindent{\bf Step 3}

\bigskip

We verify, consequential on Step 2, that Proposition \ref{prop_prox2} is available for $v_k$. I.e., there exists a function $h\in\mathcal{C}^{2,{\gamma}}_{loc}(B_1)$ such that
\[
	\left\|v_k\,-\,h\right\|_{L^\infty(B_{1/2})}\,\leq\,\delta.
\]
Set 
\[	
	\overline{P}_k(x)\,:=\,h(0)\,+Dh(0)\cdot x\,+\,\frac{1}{2}x^TD^2h(0)x
\]
and apply the triangular inequality to find that
\begin{equation}\label{eq_triangular}
	\sup_{B_\rho}\,|v_k(x)\,-\,\overline{P}_k(x)|\,\leq\,\delta\,+\,C\rho^{2+\gamma}.
\end{equation}
Now, we make the following (universal) choices for both $\delta$ and $\rho$:
\begin{equation}\label{eq_univchoicesc1ll}
	\delta\,:=\,\frac{\rho^{2}}{2}\;\;\;\;\;\;\;\;\;\;\mbox{and}\;\;\;\;\;\;\;\;\;\;\rho\,:=\,\left(\frac{1}{2C}\right)^\frac{1}{{\gamma}}.
\end{equation}
By combining \eqref{eq_triangular} with \eqref{eq_univchoicesc1ll}, we obtain
\[
	\sup_{B_\rho}\,|v_k(x)\,-\,\overline{P}_k)|\,\leq\,\rho^{2},
\]
which amounts to
\begin{equation}\label{eq_condp1}
	\sup_{B_{\rho^{k+1}}}\,\left|u(x)\,-\,[P_k(x)\,+\,\rho^{k}\overline{P}_k(x)]\right|\,\leq\,\rho^{2(k+1)}.
\end{equation}
Finally, set
\[
	P_{k+1}(x)\,:=\,P_k(x)\,+\,\rho^{k}\overline{P}_k(x\,-\,x_0).
\]

\bigskip

\noindent{\bf Step 4}

\bigskip

The definition of $P_{k+1}$ together with \eqref{eq_condp1} produces the $k+1$-th step for \eqref{eq_condc1ll2}. In addition, $C_{k+1}=C_k+D^2h(0)$; therefore,
\[
	\inf_{\beta\in\mathcal{B}}\,\tr\left(\overline{A}_\beta(0) \,C_{k+1}\right)\,=\,\left\langle f\right\rangle,
\]
which verifies the case $n=k+1$ for \eqref{eq_condc1ll1}. To conclude the proof, we note that 
\[
	a_{k+1}=a_k+\rho^{2k}h(0),\;\;\;\;\;\;\;\bb_{k+1}=\bb_k+\rho^kDh(0)\;\;\;\;\;\mbox{and}\;\;\;\;\;C_{k+1}=C_k+D^2h(0);
\]
hence, the $(k+1)$-th step for \eqref{eq_condc1ll3} amounts to ensure that
\[
	|h(0)|\,+\,|Dh(0)|\,+\,|D^2h(0)|\,\leq\,C.
\]
However, this inequality follows from the regularity of $h$.
\end{proof}

\begin{proof}[Proof of Theorem \ref{teo_main2}]
Without loss of generality, we prove the result at the origin. The proof is consequential on the following observation. From \eqref{eq_condc1ll3} we conclude that $(a_n)_{n\in\mathbb{N}}$ and $(\bb_n)_{n\in\mathbb{N}}$ are convergent sequences. Moreover,
\[
	|a_n\,-\,u(0)|\,\leq\,C\rho^{2n}\;\;\;\;\;\;\;\;\;\;\mbox{and}\;\;\;\;\;\;\;\;\;\;|\bb_n\,-\,Du(x)|\,\leq\,C\rho^n;
\]
hence, $a_n\to u(x_0)$ whereas $\bb_n\to Du(x_0)$. As regards, the sequence $(C_n)_{n\in\mathbb{N}}$, we claim that
\[
	|C_n|\,\leq\,nC;
\]
in fact,
\[
	|C_1|\,=\,|C_1\,-\,C_0|\,\leq\,C.
\]
If the case $n=k$ has been verified, we have
\[
	|C_{k+1}|\,\leq\,C\,+\,|C_k|\,\leq\,C\,+\,kC\,=\,(k+1)C.
\]

Let $r\in(0,1/2)$ be given. Choose $k\in\mathbb{N}$ so that $\rho^{k+1}< r\leq \rho^k$. The previous computations yield
\begin{align*}
	\sup_{B_r(0)}|u(x)-[u(0)+Du(0)\cdot x]|&\leq \rho^{2k}+|u(0)-a_k|+\rho|Du(0)-\bb_k|\\&\quad+\rho^{2k}|C_k|\\& \leq Cr^2\ln r^{-1},
\end{align*}
producing the result at the origin. A change of variables concludes the proof.
\end{proof}

\begin{Remark}\label{remdubmo}
Minor modifications to the ideas in Remark \ref{rem_bmo} build upon Theorem \ref{teo_main2} to yield $p-\bmo$ regularity for $Du$. The key aspect in this setting is to produce a sequence of approximating \textit{affine} functions, instead of approximating polynomials. Refer to \cite{teixeirauniversal} for further details.
\end{Remark}

\section{Estimates in $\mathcal{C}^{2,\alpha}_{loc}(B_1)$}\label{isaacsc2alpha}

In this section, we detail the proof of Theorem \ref{teo_main3}. We start with a proposition.

\begin{Proposition}\label{prop_polinomio}
Let $u\in\mathcal{C}(B_1)$ be a viscosity solution to \eqref{eq_isaacs}. Suppose A\ref{assump_ellipticity} and A\ref{assump_oscc2alpha} are in force. There exists a sequence of polynomials $(P_n)_{n\in\mathbb{N}}$, of the form
\[
	P_n(x)\,:=\,a_n\,+\,\bb_n\cdot x\,+\,\frac{1}{2}x^TC_nx,
\]
with $P_0\equiv P_{-1}\equiv 0$, and a number $\rho\in(0,1/2)$ satisfying
\begin{equation}\label{eq_condc2alpha1}
	\inf_{\beta\in\mathcal{B}}\,\left[-\tr\left(\overline{A}_\beta C_n\right)\right]\,=\,0,
\end{equation}
\begin{equation}\label{eq_condc2alpha2}
	\left\|u\,-\,P_n\right\|_{L^\infty(B_{\rho^n})}\,\leq\,\rho^{n(2+\gamma)}
\end{equation}
and
\begin{equation}\label{eq_condc2alpha3}
	|a_n-a_{n-1}|_+\rho^{n-1}|\bb_n-\bb_{n-1}|+\rho^{2(n-1)}\left\|C_n-C_{n-1}\right\|\,\leq\,C\rho^{(n-1)(2+\gamma)},
\end{equation}
for every $n\geq 0$.
\end{Proposition}
\begin{proof}
We present the proof in four steps. As before, we argue by induction in $n\geq 0$. The case $n=0$ is obvious. Suppose we have established the case $n=k$. Next we study the case $n=k+1$.

\bigskip

\noindent{\bf Step 1}

\bigskip

Consider the function $v_k:B_1\to\mathbb{R}$ defined by
\[
	v_k(x)\,:=\,\frac{(u\,-\,P_k)(\rho^kx)}{\rho^{k(2+\gamma)}}.
\]
Notice that $v_k$ solves
\begin{equation}\label{scaledc2alpha}
	\frac{1}{\rho^{\gamma k}}\sup_{\alpha\in\mathcal{A}}\inf_{\beta\in\mathcal{B}}\,\left[-\tr\left(A_{\alpha,\beta}(\rho^kx)(\rho^{\gamma k}D^2v_k(x)\,+\,C_k)\right)\right]\,=\,f_k(x)
	\end{equation}
with
\[
	f_k(x)\,:=\,\frac{f(\rho^kx)}{\rho^{\gamma k}}.
\]

\bigskip

\noindent{\bf Step 2}

\bigskip

In what follows, we approximate $v_k$ by a suitable function. Set $M_k:=\rho^{\gamma k}M+C_k$ and compute
\begin{align}
\nonumber&\left|	\sup_{\alpha\in\mathcal{A}}\inf_{\beta\in\mathcal{B}}\left[-\tr\left(A_{\alpha,\beta}(\rho^kx)M_k\right)\right]-\inf_{\beta\in\mathcal{B}}\left[-\tr\left(\overline{A}_\beta M_k\right)\right]\right|\\
\nonumber&\qquad	\leq \sup_{\alpha\in\mathcal{A}}\sup_{\beta\in\mathcal{B}}\left|\tr\left[\left(A_{\alpha,\beta}(\rho^kx)-\overline{A}_\beta\right)M_k\right]\right|\\
\nonumber&\qquad    \leq \sup_{x\in B_1}\sup_{\alpha\in\mathcal{A}}\sup_{\beta\in\mathcal{B}}\left|A_{\alpha,\beta}(\rho^kx)-\overline{A}_\beta\right|\left\|\rho^{\gamma k}M+C_k\right\|\\
\nonumber&\qquad 	\leq C(d)\epsilon_3\rho^{\gamma k}\left\|\rho^{\gamma k}M+C_k\right\|\\
\nonumber&\qquad \leq C(d)\epsilon_3\rho^{\gamma k}\left(\left\|M\right\|\,+\,\left\|C_k\right\|\right),
\end{align}
where the third inequality follows from A\ref{assump_oscc2alpha}. Because \eqref{eq_condc2alpha3} has been checked for $n=k$, we have
\[
	\left\|C_k\right\|\,\leq\,\frac{C(1-\rho^{k-1})}{1-\rho}\,\leq\,\frac{C}{1-\rho}\,\leq\, \widetilde{C}.
\]
Hence,
\begin{align}\label{eq_convuniflocc2alpha}
\nonumber	\frac{1}{\rho^{\gamma k}}&\left|	\sup_{\alpha\in\mathcal{A}}\inf_{\beta\in\mathcal{B}}\left[-\tr\left(A_{\alpha,\beta}(\rho^kx)M_k\right)\right]-\inf_{\beta\in\mathcal{B}}\left[-\tr\left(\overline{A}_\beta M_k\right)\right]\right|\\
		&\qquad\;\;\;\;\;\leq C(d)\widetilde{C}\epsilon_3\left(1\,+\left\|M\right\|\right).
\end{align}
In addition, A\ref{assump_oscc2alpha} yields
\begin{equation}\label{eq_smallf}
	\left\|f_k\right\|^p_{L^p(B_1)}\,=\,\frac{1}{\rho^{\gamma k}}\int_{B_1}\left|f(\rho^kx)\right|^pdx\,=\,\frac{1}{\rho^{\gamma k}}\intav{B_{\rho^{\gamma k}}}\left|f(y)\right|^pdy\,\leq\,\epsilon_3^p.
\end{equation}

In the sequel, we gather \eqref{eq_convuniflocc2alpha} and \eqref{eq_smallf} with standard stability results for viscosity solutions. Then, we conclude that for every $\delta>0$ there is a choice of $\epsilon_3=\epsilon_3(\delta)$ that ensures the existence of $h\in\mathcal{C}(B_{8/9})$ solving
\begin{equation}\label{eq_approxc2alpha}
	\begin{cases}
		\frac{1}{\rho^{\gamma k}}\inf_{\beta\in\mathcal{B}} \left[-\tr\left(\overline{A}_\beta(\rho^{\gamma k}D^2h(x)+C_k)\right)\right]\,=\,0&\;\;\;\;\;\mbox{in}\;\;\;\;\;B_{8/9}\\
		h\,=\,v_k&\;\;\;\;\;\mbox{on}\;\;\;\;\;\partial B_{8/9},
	\end{cases}
\end{equation}
with
\begin{equation}\label{eq_vkhc2alpha}
	\left\|v_k\,-\,h\right\|_{L^\infty(B_{8/9})}\,\leq\,\delta.
\end{equation}
The (universal) choice of $\delta$ in the next step of this proof determines $\epsilon_3$ in A\ref{assump_oscc2alpha}. 

Following, we examine the regularity of $h$. We notice that 
\begin{equation}\label{eq_auxc2alphax}
	\inf_{\beta\in\mathcal{B}}\left[-\tr\left(\overline{A}_\beta \frac{C_k}{\rho^{\gamma k}}\right)\right]\,=\,\frac{1}{\rho^{\gamma k}}\inf_{\beta\in\mathcal{B}}\left[-\tr\left(\overline{A}_\beta C_k\right)\right]\,=\,0,
\end{equation}
where the second equality follows from \eqref{eq_condc2alpha1} together with the induction hypothesis. Moreover,
\begin{align*}
	0\,=\,\frac{1}{\rho^{\gamma k}}\inf_{\beta\in\mathcal{B}}\left[-\tr\left(\overline{A}_\beta\left(\rho^{\gamma k}D^2h+C_k\right)\right)\right]\,=\,\inf_{\beta\in\mathcal{B}}\left[-\tr\left(\overline{A}_\beta\left(D^2h+\frac{C_k}{\rho^{\gamma k}}\right)\right)\right].
\end{align*}
The Evans-Krylov theory implies that
\[
	\inf_{\beta\in\mathcal{B}}\left[-\tr\left(\overline{A}_\beta D^2v\right)\right]\,=\,0\;\;\;\;\;\mbox{in}\;\;\;\;\;B_{8/9}
\]
has $\mathcal{C}^{2,\overline{\gamma}}$-estimates with
\[
	\left\|v\right\|_{\mathcal{C}^{2,\overline{\gamma}}(B_{1/2})}\,\leq\, C,
\]
where $\overline{\gamma}\in(0,1)$ and $C>0$ are universal constants. Combining this fact with \eqref{eq_auxc2alphax}, we conclude that $h\in\mathcal{C}^{2,\overline{\gamma}}_{loc}(B_1)$ with 
\[
	\left\|h\right\|_{\mathcal{C}^{2,\overline{\gamma}}(B_{1/2})}\,\leq\,C.
\]
Then,
\begin{equation}\label{eq_hpolc2alpha}
	\left\|h-\left[h(0)+Dh(0)\cdot x+\frac{1}{2}x^TD^2h(0)x\right]\right\|_{L^\infty(B_\rho)}\,\leq\,C\rho^{2+\overline{\gamma}}.
\end{equation}

\bigskip

\noindent{\bf Step 3}

\bigskip

Now we combine \eqref{eq_vkhc2alpha} with \eqref{eq_hpolc2alpha} and use the triangular inequality to show that
\begin{align*}
	\left\|v_k-\left[h(0)+Dh(0)\cdot x+\frac{1}{2}x^TD^2h(0)x\right]\right\|_{L^\infty(B_\rho)}\,\leq\,\delta\,+\,C\rho^{2+\overline{\gamma}}.
\end{align*}
By setting
\[
	\delta\,:=\,\frac{\rho^{2+\gamma}}{2},\;\;\;\;\;\;\;\;\;\;\;\;\;\;\;\;\;\;\;\;\rho\,:=\,\left(\frac{1}{2C}\right)^\frac{1}{\overline{\gamma}-\gamma}
\]
and
\[
	\overline{P}_k(x)\,:=\,h(0)+Dh(0)\cdot x+\frac{1}{2}x^TD^2h(0)x,
\]
we obtain
\begin{equation}\label{eq_vkmenospol}
	\left\|v_k\,-\,\overline{P}_k\right\|_{L^\infty(B_\rho)}\,\leq\,\rho^{2+\gamma}.
\end{equation}

\bigskip

\noindent{\bf Step 4}

\bigskip

In what follows, we define the polynomial $P_{k+1}$ and conclude the proof. In light of the definition of $v_k$, \eqref{eq_vkmenospol} leads to
\[
	\sup_{x\in B_{\rho^{k+1}}}\left|u(x)\,-\,P_{k}(x)\,-\,\rho^{k(2+\gamma)}\overline{P}_k(\rho^{-k}x)\right|\,\leq\,\rho^{(k+1)(2+\gamma)}.
\]
Therefore, by setting $P_{k+1}(x)\,:=\,P_k(x)+\rho^{k(2+\gamma)}\overline{P}_k(\rho^{-k}x)$, we get
\[
	\left\|u\,-\,P_{k+1}\right\|_{L^\infty(B_{\rho^{k+1}})}\,\leq\,\rho^{(k+1)(2+\gamma)}.
\]
This verifies \eqref{eq_condc2alpha2}. As regards \eqref{eq_condc2alpha1}, notice that $C_{k+1} = C_k+\rho^{\gamma k}D^2h(0)$; hence, \eqref{eq_approxc2alpha} implies
\[
	\inf_{\beta\in\mathcal{B}}\,\left[-\tr\left(\overline{A}_\beta C_{k+1}\right)\right]\,=\,0.
\]
To verify \eqref{eq_condc2alpha3}, observe that 
\[
	\left|a_{k+1}\,-\,a_k\right|\,\leq\,\rho^{k(2+\alpha)}h(0),
\]
\[
	\rho^k\left|\bb_{k+1}\,-\,\bb_k\right|\,\leq\,\rho^{k(2+\alpha)}Dh(0),
\]
and
\[
	\rho^{2k}\left\|C_{k+1}\,-\,C_k\right\|\,\leq\,\rho^{k(2+\alpha)}D^2h(0);
\]
finally, use the $\mathcal{C}^{2,\overline{\gamma}}$-estimates available for $h$. This finishes the proof.

\end{proof}

We close this section with the proof of Theorem \ref{teo_main3}.

\begin{proof}[Proof of Theorem \ref{teo_main3}]
We first prove Assertion $1.$; in view of Proposition \ref{prop_polinomio}, there exists a polynomial $P^*$ such that $P_n\to P^*$, uniformly in $B_1$. The regularity of the approximate function $h$ ensures that
\[
	|DP^*(0)|\,+\,\left\|D^2P^*(0)\right\|\,\leq\,C.
\]
Moreover, for every $n\in\mathbb{N}$ we have
\[
	\left\|u\,-\,P^*\right\|_{L^\infty(B_{\rho^n})}\,\leq\, C\rho^{n(2+\gamma)};
\]
see \cite[Chapter 8; p.76]{ccbook}. For Assertion $2.$, we notice the result follows from a change of variables argument. This concludes the proof of the theorem.
\end{proof}






\bibliography{bib_2016}
\bibliographystyle{plain}

\bigskip


\noindent\textsc{Edgard A. Pimentel}\\
Department of Mathematics\\
Pontifical Catholic University of Rio de Janeiro -- PUC-Rio\\
22451-900, G\'avea, Rio de Janeiro-RJ, Brazil\\
\noindent\texttt{pimentel@puc-rio.br}

\end{document}